\documentclass[11pt]{amsart}
\usepackage[all]{xy}
\usepackage{amssymb,mathrsfs,euscript,enumitem,amsmath,bbold}

\usepackage[in]{fullpage}
\usepackage{color}
\definecolor{Chocolat}{rgb}{0.36, 0.2, 0.09}
\definecolor{BleuTresFonce}{rgb}{0.215, 0.215, 0.36}
\definecolor{EgyptianBlue}{rgb}{0.06, 0.2, 0.65}

\usepackage[colorlinks,final,hyperindex]{hyperref}
\hypersetup{citecolor=BleuTresFonce, urlcolor=EgyptianBlue, linkcolor=Chocolat}

\usepackage{fourier}

\DeclareMathAlphabet{\mathbbold}{U}{bbold}{m}{n}

\def\k{\mathbbold{k}}

\DeclareMathOperator{\Tor}{Tor}

\newtheorem*{mtheorem}{Theorem}
\newtheorem{theorem}{Theorem}

\newtheorem*{conjecture}{Conjecture}

\theoremstyle{definition}

\newtheorem{remark}[theorem]{Remark}

\begin{document}

\title{Anick resolution and Koszul algebras of finite global dimension}

\author{Vladimir Dotsenko}
\address{School of Mathematics, Trinity College, Dublin 2, Ireland}
\email{vdots@maths.tcd.ie}

\author{Soutrik Roy Chowdhury}
\address{ITI Road, Jyotinagar, Siliguri 734001, India}
\email{roychowdhurysoutrik@gmail.com}


\begin{abstract}
We show how to study a certain associative algebra recently discovered by Iyudu and Shkarin using the Anick resolution. This algebra is a counterexample to the conjecture of Positselski on  Koszul algebras of finite global dimension. 
\end{abstract}

\maketitle

In the definitive treatment of quadratic algebras, the monograph \cite{PP}, the following conjecture of Positselski is recorded:

\begin{conjecture}[Conjecture 2 in {\cite[Chapter~7]{PP}}]
Any Koszul algebra $A$ of finite global homological dimension $d$ has $\dim A_1\ge d$. Dually, for a Koszul algebra $B$ with $B_{d+1}=0$ and $B_d\ne 0$ one has $\dim B_1\ge d$. 
\end{conjecture}

In a recent preprint of Iyudu and Shkarin \cite{IS} where a classification result for Hilbert series of Koszul algebras with three generators and three relations is proposed, one of the algebras from the classification theorem provides a counterexample to this conjecture. Namely, the following result holds.

\begin{mtheorem}
Consider the following algebra with three generators and three quadratic relations:
 \[
A=\k\langle x,y,z\mid x^2+yx, xz,zy\rangle . 
 \]
The algebra $A$ is Koszul. Its global dimension is equal to $4$.  Hence, the conjecture above is false. 
\end{mtheorem}

The proof of Koszulness of the algebra $A$ in \cite[Prop. 4.2]{IS}, contains a typo of a sort: what is claimed to be the reduced Gr\"obner basis for $A$ is not really a Gr\"obner basis. The strategy of \cite{IS} appears to be valid if one uses the correct reduced Gr\"obner basis; however, we would like to present a completely different proof of this result which illustrates the power of the Anick resolution for associative algebras~\cite{Anick}. 

\begin{proof}

Throughout this proof, we use the degree-lexicographic ordering of variables with $x>y>z$. It turns out that the reduced Gr\"obner basis of $A$ for this ordering, though infinite, has a very simple and pleasant description. Namely, it consists of the following polynomials:
 \[
xy^kx+y^{k+1}x (k\ge 0), xz, zy .
 \]
To establish that, we use the Diamond Lemma \cite{Bergman}. First of all, we note that the ideal generated by the relation $x^2+yx$ alone already contains all the elements $xy^kx+y^{k+1}x$ above, and that the overlaps of their leading terms give S-polynomials that can be reduced to zero modulo these elements; it is a computation essentially identical to the one of \cite[Ex. 3.6.2]{RC-MSc}. Finally, we note that the S-polynomial of $xy^kx+y^{k+1}x$ and $xz$ is 
 \[
(xy^kx+y^{k+1}x)z-xy^k(xz)=y^{k+1}xz , 
 \]
which is divisible by $xz$ and hence can be reduced to zero, and the S-polynomial of $xz$ and $zy$ is actually equal to zero, as those relations are monomial. 

Furthermore, it turns out that the Anick resolution \cite{Anick} for the Gr\"obner basis we computed is surprisingly manageable. We shall use the description of that resolution due to Ufnarovski \cite[Sec.~3.6]{Ufn}, which we recall here for the sake of completeness. That definition recursively defines Anick chains and their tails for any algebra $A=\k\langle X\mid R\rangle$, where the set of defining relations is a Gr\"obner basis:
\begin{itemize}
\item[-] elements of $X$ are the only $0$-chains, each of them coincides with its tail;
\item[-] each $q$-chain is word $c$ in the alphabet $X$ which is the concatenation $c't$, where $t$ is the tail of $c$, and $c'$ is a $(q-1)$-chain;
\item[-] if we denote by~$t'$ the tail of $c'$ in the above decomposition, there exists a ``factorization'' $t'=m_1m_2$  with $m_2t$ is the leading term of an element of $R$, 
and there are no other divisors of $t't$ that are leading terms of $R$.
\end{itemize}
In our case, a direct inspection of the leading terms $xy^kx$, $xz$, $zy$ of the reduced Gr\"obner basis we found instantly produces the list of all Anick chains as follows. 
\begin{itemize}
\item the set of $0$-chains $C_0$ consists of $x,y,z$;
\item the set of $1$-chains $C_1$ consists of $xy^kx (k\ge 0)$, $xz$, $zy$;
\item for $n\ge 2$, the set of $n$-chains $C_n$ consists of $xy^{k_1}x\cdots xy^{k_n}x$ $(k_1, \ldots, k_n\ge 0)$, $xy^{k_1}x\cdots xy^{k_{n-1}}xz$ $(k_1, \ldots, k_{n-1}\ge 0)$, $xy^{k_1}x\cdots xy^{k_{n-2}}xzy$ $(k_1, \ldots, k_{n-2}\ge 0)$.
\end{itemize}
It is known \cite{Anick,Ufn} that there exists a resolution of the (right) augmentation $A$-module $\k$ by free right $A$-modules of the form 
 \[
\ldots\to \k C_n\otimes A\to \k C_{n-1}\otimes A\to\ldots \to \k C_0\otimes A\to A\to 0. 
 \]
The general recipe for the computation of the differential of that resolution in our case gives the following results. First, the formulas for $d_0$ and $d_1$ are standard:
\begin{gather*}
d_0(x\otimes 1)=x,
d_0(y\otimes 1)=y,
d_0(z\otimes 1)=z,\\
d_1(xy^kx\otimes 1)=x\otimes y^kx+y\otimes y^kx,
d_1(xz\otimes 1)=x\otimes z, 
d_1(zy\otimes 1)=z\otimes y. 
\end{gather*}
Computing $d_n$ with $n\ge 2$ is very similar to the computation of \cite[Prop.~4.3.1]{RC-MSc}. Namely, the following formulas hold:
\begin{multline*}
d_n(xy^{k_1}x\cdots xy^{k_n}x\otimes 1)=\\ 
xy^{k_1}x\cdots xy^{k_{n-1}}x\otimes y^{k_n}x+
\sum_{i=1}^{n-1}(-1)^{n-1-i}xy^{k_1}x\cdots xy^{k_{i-1}}xy^{k_i+k_{i+1}+1}x y^{k_{i+2}}x\cdots xy^{k_n}x\otimes 1 ,
\end{multline*}
\begin{multline*}
d_n(xy^{k_1}x\cdots xy^{k_{n-1}}xz\otimes 1)=\\ 
xy^{k_1}x\cdots xy^{k_{n-1}}x\otimes z+
\sum_{i=1}^{n-2}(-1)^{n-2-i}xy^{k_1}x\cdots xy^{k_{i-1}}xy^{k_i+k_{i+1}+1}x y^{k_{i+2}}x\cdots xy^{k_{n-1}}xz\otimes 1 ,
\end{multline*}
\begin{multline*}
d_n(xy^{k_1}x\cdots xy^{k_{n-2}}xzy\otimes 1)=\\ 
xy^{k_1}x\cdots xy^{k_{n-1}}xz\otimes y+
\sum_{i=1}^{n-3}(-1)^{n-3-i}xy^{k_1}x\cdots xy^{k_{i-1}}xy^{k_i+k_{i+1}+1}x y^{k_{i+2}}x\cdots xy^{k_{n-2}}xzy\otimes 1 .
\end{multline*}
To compute the bar homology of $A$, or equivalently $\Tor^A_\bullet(\k,\k)$, we should compute the homology of the complex 
 \[
\Bigl\{\left(\k C_n\otimes A\right)\otimes_A \k,\, \bar{d}_n\Bigr\} ,
 \]
where the induced differential $\bar{d}_n=d_n\otimes 1$ only retains the terms that are not annihilated by the augmentation of $A$. In other words, under the identification $\left(\k C_n\otimes A\right)\otimes_A \k\cong \k C_n$ we have 
\begin{gather*}
\bar{d}_0(x)=0,\quad 
\bar{d}_0(y)=0,\quad 
\bar{d}_0(z)=0,\\
\bar{d}_1(xy^kx)=0,\quad
\bar{d}_1(xz)=0, \quad 
\bar{d}_1(zy)=0, \\
\bar{d}_n(xy^{k_1}x\cdots xy^{k_n}x)=
\sum_{i=1}^{n-1}(-1)^{n-1-i}xy^{k_1}x\cdots xy^{k_{i-1}}xy^{k_i+k_{i+1}+1}x y^{k_{i+2}}x\cdots xy^{k_n}x ,\\
\bar{d}_n(xy^{k_1}x\cdots xy^{k_{n-1}}xz)=
\sum_{i=1}^{n-2}(-1)^{n-2-i}xy^{k_1}x\cdots xy^{k_{i-1}}xy^{k_i+k_{i+1}+1}x y^{k_{i+2}}x\cdots xy^{k_{n-1}}xz ,\\
\bar{d}_n(xy^{k_1}x\cdots xy^{k_{n-2}}xzy)= 
\sum_{i=1}^{n-3}(-1)^{n-3-i}xy^{k_1}x\cdots xy^{k_{i-1}}xy^{k_i+k_{i+1}+1}x y^{k_{i+2}}x\cdots xy^{k_{n-2}}xzy .
\end{gather*}
These formulas show that the chain complex we are dealing with decomposes as a direct sum of the subcomplex
 \[
0\to \k\{x^2zy\}\to \k\{x^2z, xzy\}\to \k\{x^2, xz, zy\} \to\k\{x,y,z\}\to\k\to 0
 \]
which has zero differential, the subcomplexes $U_\bullet^{(p)}$, $p\ge1$, spanned by the elements $xy^{k_1}x\cdots xy^{k_n}x$, $n\ge 1$, the subcomplexes $V_\bullet^{(p)}$, $p\ge 1$, spanned by the elements $xy^{k_1}x\cdots xy^{k_{n-1}}xz$, $n\ge 2$, and the subcomplexes $W_\bullet^{(p)}$, $p\ge 1$, spanned by the elements $xy^{k_1}x\cdots xy^{k_{n-2}}xzy$, $n\ge 3$; in each of these cases, the superscript $(p)$ imposes the constraint $\sum\limits_{i} (k_i+1)=p+1$. We note that $V_\bullet^{(p)}\cong U_{\bullet+1}^{(p)}$ and $W_\bullet^{(p)}\cong U_{\bullet+2}^{(p)}$, so if we prove that each subcomplex $U_\bullet^{(p)}$, $p\ge1$, is acyclic, this will prove that all other complexes we consider are acyclic as well.

Note that the subcomplex $U_\bullet^{(p)}$ appears when computing the bar homology of the algebra $A'=\k\langle x,y\mid x^2=yx\rangle$ via the Anick resolution. Indeed, after tensoring the Anick resolution of the augmentation module for $A'$ with the augmentation module, the resulting complex decomposes as a direct sum of the complex
 \[
0\to \k\{x^2\} \to\k\{x,y\}\to\k\to 0
 \]
with zero differential and complexes $U_\bullet^{(p)}$ with $p\ge 1$. At the same time, one can note that for the algebra $A'$, the degree-lexicographic order with $x<y$ gives a quadratic Gr\"obner basis, yielding the complex 
 \[
0\to \k\{yx\} \to\k\{x,y\}\to\k\to 0
 \]
with zero differential which represents the bar homology of $A'$. This implies that $A'$ is Koszul, and that its bar homology for the ordering we are interested in must be represented by the classes of the Anick chains $1,x,y,x^2$. Therefore, all complexes $U_\bullet^{(p)}$ with $p\ge 1$ must be acyclic. 

It follows that the subcomplex 
 \[
0\to \k\{x^2zy\}\to \k\{x^2z, xzy\}\to \k\{x^2, xz, zy\} \to\k\{x,y,z\}\to\k\to 0
 \]
we noted above represents the bar homology of $A$. Thus, the bar homology is concentrated on the diagonal, so $A$ is Koszul. Moreover, its Koszul dual coalgebra vanishes in degrees higher than four and is nonzero in degree four, which by the classical results on the Koszul duality \cite{BGS88} proves the claim on the global dimension, therefore completing the proof.
\end{proof}
 
\begin{remark}
The fact that the complexes $U_\bullet^{(p)}$ are acyclic can be also established by identifying them with graded pieces of the Shafarevich complex \cite{GS64,G88} for the subset $\{t\}$ of $\k\langle t\rangle$, but we chose the argument above to emphasize how one can use the Anick resolutions for two different orderings simultaneously. 
\end{remark}

\subsection*{Acknowledgements} 
We are grateful to Leonid Positselski for drawing our attention to the fact that the list of algebras in \cite{IS} contains a counterexample to his conjecture. Ivan Yudin informed us that it is also possible to establish the claim on Koszulness by a careful examination of the two Anick resolutions for the orderings $x>y>z$ and $y>x>z$ for the original algebra; his argument utilises Ufnarovski's graph for generating Anick chains \cite{Ufn,Ufn1}; we are very thankful for that remark.

\bibliographystyle{amsplain}

\begin{thebibliography}{10}
\bibitem{Anick} David~J. Anick, \emph{On the homology of associative algebras}, Trans. Amer. Math. Soc. \textbf{296} (1986), no.~2, 641--659.
\bibitem{BGS88} Alexander Beilinson, Victor Ginsburg, and Vadim Schechtman, \emph{Koszul duality}, J. Geom. Phys. \textbf{5} (1988), no.~3, 317--350.
\bibitem{Bergman} George Bergman, \emph{The diamond lemma for ring theory}, Adv. in Math. \textbf{29} (1978), no. 2, 178--218. 
\bibitem{GS64} Evgeny S.~Golod and Igor R.~Shafarevich, \emph{On the class field tower}, Izv. Akad. Nauk SSSR Ser. Mat., \textbf{28} (1964), no.~2, 261--272.
\bibitem{G88}  Evgeny S. Golod, \emph{Standard bases and homology}, Lecture Notes in Math. \textbf{1352} (1988), 88--95.
\bibitem{IS} Nataliya Iyudu and Stanislav Shkarin, \emph{One question from the Polishchuk and Positselski book on Quadratic algebras}, Preprint IHES/M/16/16, \url{http://preprints.ihes.fr/2016/M/M-16-16.pdf}, May 2016. 
\bibitem{PP} Alexander Polishchuk and Leonid Positselski, \emph{Quadratic algebras}, University Lecture Series, \textbf{37}, AMS, Providence, RI, 2005.
\bibitem{RC-MSc} Soutrik Roy Chowdhury, \emph{Gr\"obner bases: connecting linear algebra with homological and homotopical algebra}, ArXiv preprint 1510.01542, 2015. 
\bibitem{Ufn} Victor~A. Ufnarovski, \emph{Combinatorial and asymptotic methods in algebra}, Algebra, VI, Encyclopaedia Math. Sci., vol.~57, Springer, Berlin, 1995,
  pp.~1--196.
\bibitem{Ufn1} Victor~A. Ufnarovski, \emph{On the use of graphs for computing a basis, growth and Hilbert series of associative algebras}, Mathematics of the USSR -- Sbornik, \textbf{68} (1991), no.~2, 417--428.
\end{thebibliography}
\providecommand{\bysame}{\leavevmode\hbox to3em{\hrulefill}\thinspace}

\end{document}